\documentclass[11pt]{article}

\usepackage[dvips]{graphics}
\usepackage{amsfonts}
\usepackage{amssymb}
\usepackage{amsthm}

\usepackage[dvips]{graphics}

\theoremstyle{plain}
\newtheorem{theorem}{Theorem}[section]
\newtheorem{lemma}[theorem]{Lemma}
\newtheorem{proposition}[theorem]{Proposition}
\newtheorem{corollary}[theorem]{Corollary}
\theoremstyle{definition}
\newtheorem{definition}[theorem]{Definition}

\newtheorem{example}[theorem]{Example}

\newtheorem{notation}[theorem]{Notation}

\newcommand{\todo}[1]{\vspace{5 mm}\par \noindent
\marginpar{\textsc{ToDo}}
\framebox{\begin{minipage}[c]{0.95 \textwidth}
\tt #1 \end{minipage}}\vspace{5 mm}\par}

\renewcommand{\todo}[1]{}

\newcommand{\idiot}[1]{\vspace{5 mm}\par \noindent
\framebox{\begin{minipage}[c]{0.95 \textwidth}
\tt #1 \end{minipage}}\vspace{5 mm}\par}

\renewcommand{\idiot}[1]{}

\newcommand{\xs}{x_1,\ldots,x_n}                
\newcommand{\vs}{v_1,\ldots,v_n}                
\newcommand{\Ms}{M_1,\ldots,M_q}                
\newcommand{\Fs}{F_1,\ldots,F_q}                
\newcommand{\Gs}{G_1,\ldots,G_p}                
\newcommand{\As}{A_1,\ldots,A_t}                
\newcommand{\Bs}{B_1,\ldots,B_s}                

\newcommand{\dimn}{ {\rm{dim}} \ }              
\newcommand{\F}{{\mathcal{F}}}                  
\newcommand{\N}{{\mathcal{N}}}                  
\newcommand{\D}{\Delta}                         
\newcommand{\DN}{\Delta_N}                      
\newcommand{\al}{\alpha}                         
\newcommand{\df}{\delta_{\mathcal{F}}}          
\newcommand{\dn}{\delta_{\mathcal{N}}}          
              
\newcommand{\tuple}[1]{\langle #1 \rangle}      
\newcommand{\rmv}[1]{\setminus \langle #1\rangle} 

\newcommand{\st}{\ | \ }
\newcommand{\DM}{\Delta_M}                      
  
\newcommand{\alx}[1]{{#1{}^{\vee}}}               
\newcommand{\IM}{\alx{I}}                       
\newcommand{\IMi}{{\alx{I}}_{[i]}}                
\newcommand{\DMi}{\Delta_{M,\ [i]}}             
\newcommand{\DNi}{\Delta_{N,\ i}}               
\newcommand{\cmp}[1]{{#1}^c}                    
\newcommand{\tor}[2]{{\rm Tor}_{#1}^{#2}\ }     

\title{\sc Simplicial Trees are Sequentially Cohen-Macaulay} \author{
Sara Faridi\thanks{Universit\'{e} du Qu\'{e}bec \`{a} Montr\'{e}al,
Laboratoire de combinatoire et d'informatique math\'{e}matique, Case
postale 8888, succursale Centre-Ville, Montr\'{e}al, QC Canada H3C
3P8.  email: \emph{faridi@math.uqam.ca.} \hspace{.5in} This research
was supported by a NSERC Postdoctoral Fellowship.}}  \date{\today}

\begin{document}

\maketitle

\begin{abstract}  This paper uses dualities between facet ideal theory and 
  Stanley-Reisner theory to show that the facet ideal of a simplicial tree is
 sequentially Cohen-Macaulay. The proof involves showing that the Alexander
 dual (or the cover dual, as we call it here) of a simplicial tree is a
 componentwise linear ideal. We conclude with additional combinatorial
 properties of simplicial trees.
\end{abstract}

The main result of the this paper is that the facet ideal of a
simplicial tree is sequentially Cohen-Macaulay. Sequentially
Cohen-Macaulay modules were introduced by Stanley [S] (following the
introduction of nonpure shellability by Bj\"{o}rner and Wachs [BW]) so
that a nonpure shellable simplicial complex had a sequentially
Cohen-Macaulay Stanley-Reisner ideal. Herzog and Hibi ([HH]) then
defined the notion of a componentwise linear ideal, which extended a
criterion of Eagon and Reiner ([ER]) for Cohen-Macaulayness of an
ideal to a criterion for sequential Cohen-Macaulayness.

Simplicial trees, on the other hand, were introduced in [F1] in the
context of Rees rings, and their facet ideals were studied further in
[F2] for their Cohen-Macaulay properties, and in [Z] for their
resolutions.  The facet ideal of a given simplicial complex is a
square-free monomial ideal where every generator is the product of the
vertices of a facet of the complex. If the simplicial complex is a
tree (Definition~\ref{tree-def}), it turns out that its facet ideal
has many interesting algebraic and combinatorial properties.

Given a square-free monomial ideal, one could consider it as the facet
ideal of one simplicial complex, and the Stanley-Reisner ideal of
another. This in a sense gives two ``languages'' to study a
square-free monomial ideal. Below we provide a dictionary which makes
it easy to move from one language to the other. We use this dictionary
to translate existing criteria for Cohen-Macaulayness and sequential
Cohen-Macaulayness into the language of facet ideals, and then finally
use these criteria to show that the facet ideal of a simplicial tree
is sequentially Cohen-Macaulay (Corollary~\ref{main-result}).

There are several byproducts. An immediate one is that the facet ideal
of an unmixed simplicial tree (Definition~\ref{unmixed}) is
Cohen-Macaulay (Corollary~\ref{F2-cor}).  This is discussed at length
and proved independently in [F2], where we introduce the concept of
``grafting'' a simplicial complex. As it turns out, any unmixed tree
is grafted, and any grafted simplicial complex is Cohen-Macaulay. This
fact, in addition to proving the statement of Corollary~\ref{F2-cor},
gives the precise combinatorial structure of a Cohen-Macaulay tree.

Another outcome is that the Stanley-Reisner complex corresponding to a
Cohen-Macaulay tree is shellable. This was known in the case of graphs
([V]). In general, shellability is only a necessary condition for
Cohen-Macaulayness.

The paper is organized as follows: Section~\ref{basic-section} reviews
the basics of facet ideal theory, introducing cover complexes. In
Section~\ref{dualities-section} we discuss how facet ideal theory
relates to Stanley-Reisner theory. In Section~\ref{trees-section} we
define simplicial trees and discuss their localization. In
Section~\ref{sequential-section} we define sequentially Cohen-Macaulay
and componentwise linear ideals, and introduce a criterion for an
ideal to be sequentially Cohen-Macaulay, which we use in
Section~\ref{trees-scm-section} to prove that trees are sequentially
Cohen-Macaulay.

For the convenience of the reader, we have included a table of
notation at the end of the paper (Figure~\ref{index}).

We would like to thank J\"urgen Herzog for raising the question of
whether simplicial trees are sequentially Cohen-Macaulay, and for an
earlier reading of this manuscript.


\section{Basic Definitions}\label{basic-section}

This section is a review of the basic definitions and notations in facet
ideal theory. Much of the material here appeared in more detail in
[F1] and [F2], except for the discussion on the cover complex.

\begin{definition}[simplicial complex, facet, subcollection and more] 
  A \emph{simplicial complex} $\Delta$ over a set of vertices $V=\{
  \vs \}$ is a collection of subsets of $V$, with the property that
  $\{ v_i \} \in \Delta$ for all $i$, and if $F \in \Delta$ then all
  subsets of $F$ are also in $\Delta$ (including the empty set). An
  element of $\Delta$ is called a \emph{face} of $\Delta$, and the
  \emph{dimension} of a face $F$ of $\Delta$ is defined as $|F| -1$,
  where $|F|$ is the number of vertices of $F$.  The faces of
  dimensions 0 and 1 are called \emph{vertices} and \emph{edges},
  respectively, and $\dimn \emptyset =-1$.  The maximal faces of
  $\Delta$ under inclusion are called \emph{facets} of $\Delta$. The
  dimension of the simplicial complex $\Delta$ is the maximal
  dimension of its facets.

We denote the simplicial complex $\Delta$ with facets $\Fs$ by
$$\Delta = \tuple{\Fs}$$
and we call $\{ \Fs \}$ the \emph{facet set}
of $\Delta$. A simplicial complex with only one facet is called a
\emph{simplex}. By a \emph{subcollection} of $\Delta$ we mean a
simplicial complex whose facet set is a subset of the facet set of
$\Delta$.
\end{definition}

\begin{definition}[connected simplicial complex] A simplicial complex
  $\D=\tuple{\Fs}$ is \emph{connected} if for every pair $i,j$, $1
  \leq i < j \leq q$, there exists a sequence of facets
  $F_{t_1},\ldots,F_{t_r}$ of $\D$ such that $F_{t_1}=F_i$,
  $F_{t_r}=F_j$ and $F_{t_s} \cap F_{t_{s+1}} \neq \emptyset$ for
  $s=1,\ldots,r-1$.
\end{definition}

\begin{definition}[facet ideal, facet complex]\label{facet-ideal} Let $k$ be a 
field and $\xs$ be a set of indeterminates, and $R=k[\xs]$ be a
polynomial ring.

\begin{itemize} 
\item Let $\Delta$ be a simplicial complex over $n$ vertices labeled
  $\vs$. We define the \emph{facet ideal} of $\D$, denoted by
  $\F(\D)$, to be the ideal of $R$ generated by square-free monomials
  $x_{i_1}\ldots x_{i_s}$, where $\{v_{i_1},\ldots, v_{i_s}\}$ is a
  facet of $\Delta$.
  
\item Let $I=(\Ms)$ be an ideal in $R$, where $\Ms$ are square-free
  monomials in $\xs$ that form a minimal set of generators for $I$. We
  define the \emph{facet complex} of $I$, denoted by $\df(I)$, to be
  the simplicial complex over a set of vertices $\vs$ with facets
  $\Fs$, where for each $i$, $F_i=\{v_j \st x_j|M_i, \ 1 \leq j \leq n
  \}$.
\end{itemize}
\end{definition}

Throughout this paper we often use a letter $x$ to denote both a
vertex of $\D$ and the corresponding variable appearing in $\F(\D)$,
and $x_{i_1}\ldots x_{i_r}$ to denote a facet of $\D$ as well as a
monomial generator of $\F(\D)$.

\begin{example}\label{main-example} If $\D$ is the simplicial complex
  $\tuple{xyz,yzu,uv}$ drawn below,
\[ \includegraphics{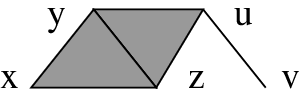} \]
  then $\F(\D)=(xyz,yuz,uv)$ is its facet ideal.
\end{example}

Facet ideals give a one-to-one correspondence between simplicial
complexes and square-free monomial ideals. 

Next we define the notion of a vertex cover. The combinatorial idea
here comes from graph theory. In algebra, it corresponds to prime
ideals lying over the facet ideal of a given simplicial complex.

\begin{definition}[vertex cover, vertex covering number, unmixed]\label{unmixed} 
  Let $\Delta$ be a simplicial complex with vertex set $V$. A
  \emph{vertex cover} for $\D$ is a subset $A$ of $V$ that intersects
  every facet of $\D$.  If $A$ is a minimal element (under inclusion)
  of the set of vertex covers of $\D$, it is called a \emph{minimal
  vertex cover}. The smallest of the cardinalities of the vertex
  covers of $\D$ is called the \emph{vertex covering number} of $\D$
  and is denoted by $\al(\D)$.
 
A simplicial complex $\Delta$ is \emph{unmixed} if all of its minimal
vertex covers have the same cardinality.
\end{definition}

\begin{example}\label{highlight-example} If $\D$ is the simplicial complex
  in Example~\ref{main-example}, then the vertex covers of $\D$ are:
$${\bf
\{x,u\},\{y,u\},\{y,v\},\{z,u\},\{z,v\}},\{x,y,u\},\{x,z,u\},\{x,y,v\},
\ldots.$$ The first five vertex covers above (highlighted in bold),
are the minimal vertex covers of $\D$.
\end{example}

In all the arguments in this paper, unless otherwise stated, $k$
denotes a field.

Given a square-free monomial ideal $I$ in a polynomial ring $k[\xs]$,
the vertices of $\df(I)$ are those variables that divide a monomial in
the generating set of $I$; this set may not include all elements of
$\{\xs\}$. The fact that some extra variables may appear in the
polynomial ring has little effect on the algebraic or combinatorial
structure of $\df(I)$. On the other hand, if $\D$ is a simplicial
complex, being able to consider the facet ideals of its subcomplexes
as ideals in the same ambient ring simplifies many of our discussions.
For this reason we make the following definition.

\begin{definition}[variable cover] Let $I$ be a square-free monomial ideal
in a polynomial ring $R=k[\xs]$. A subset $A$ of the variables
$\{\xs\}$ is called a \emph{(minimal) variable cover} of $\D=\df(I)$
(or of $I$) if $A$ is the generating set for a (minimal) prime ideal of
$R$ containing $I$.\end{definition}

If $\xs$ are all vertices of $\D=\df(I)$, then a variable cover of
$\D$ is exactly the same as a vertex cover of $\D$. In general every
variable cover of $\D$ contains a vertex cover of $\D$.  For example
for the ideal $I=(xy,xz) \subseteq k[x,y,z,u]$, $\{x,u\}$ is a
variable cover but not a vertex cover. The minimal vertex covers of
$\D$, however, are always the same as the minimal variable covers of
$\D$.

We now construct a new simplicial complex using the minimal vertex
covers of a given simplicial complex.

\begin{definition}[cover complex]\label{cover-complex} Given a simplicial 
complex $\D$, the simplicial complex $\D_M$ called the \emph{cover
  complex of} $\D$, is the simplicial complex whose facets are the
  minimal vertex covers of $\D$.
\end{definition}

\begin{example} In Example~\ref{highlight-example}, 
$\D=\tuple{xyz,yzu,uv}$ and $\D_M=\tuple{xu,yu,yv,zu,zv}.$
\end{example}

It is worth observing that $\D$ being unmixed is equivalent to $\DM$
being \emph{pure} (meaning that all facets of $\DM$ are of the same
dimension).  This fact becomes useful in our discussions below.  For
example the simplicial complex $\D$ in Example~\ref{highlight-example}
is unmixed, and $\DM$ is pure.

The following fact is known in hypergraph theory (see, for example, [B]).
We outline a proof below.

\begin{proposition}[The cover complex is a dual]\label{cover-complex-dual}
  If $\D$ is a simplicial complex, then $\D_M$ is a dual of $\D$; i.e.
  $\D_{MM}=\D.$
\end{proposition}

            \begin{proof} Suppose that $\D=\tuple{\Fs}$ and
              $\DM=\tuple{\Gs}$.  Suppose that for $i=1,\ldots,p$,
             $q_i$ is the prime ideal generated by the elements of $G_i$,
             so that we have $$I=\F(\D)= q_1 \cap \ldots \cap q_p.$$

             We first show that every facet of $\D$ is a vertex cover
             of $\DM$. Consider the facet $F_1$. Since the monomial
             $F_1 \in q_i$ for all $i$, it follows that $F_1$ contains
             at least one vertex of each of the $G_i$. This proves
             that $F_1$ is a vertex cover of $\DM$.
             
             Suppose now that $F$ is any minimal vertex cover of
             $\DM$.  Since $F$ contains a vertex of each of the $G_i$,
             it belongs to all the ideals $q_i$ (if we consider $F$ as
             a monomial), and therefore $F \in I$. So some generator
             $F_j$ of $I$ must divide $F$. This means that $F_j
             \subseteq F$, but since $F_j$ is already a vertex cover
             of $\DM$, it follows that $F=F_j$. This shows that $\Fs$
             are all the minimal vertex covers of $\DM$.
           \end{proof}


\section{Relations to Stanley-Reisner theory}\label{dualities-section}

We begin by the basic definitions from Stanley-Reisner theory. For a
detailed coverage of this topic, we refer the reader to [BH].

\begin{definition}[nonface ideal, nonface complex]\label{nonface-ideal} 
Let $k$ be a field and $\xs$ be a set of indeterminates, and
$R=k[\xs]$ be a polynomial ring.

\begin{itemize} 
\item Let $\Delta$ be a simplicial complex over $n$ vertices labeled
  $\vs$. We define the \emph{nonface ideal} or the
  \emph{Stanley-Reisner ideal} of $\Delta$, denoted by $\N(\Delta)$,
  to be the ideal of $R$ generated by square-free monomials
  $x_{i_1}\ldots x_{i_s}$, where $\{v_{i_1},\ldots, v_{i_s}\}$ is not
  a face of $\Delta$.

\item Let $I=(\Ms)$ be an ideal in $R$, where $\Ms$ are square-free
  monomials in $\xs$ that form a minimal set of generators for $I$. 
We define the \emph{nonface complex} or the \emph{Stanley-Reisner
complex} of $I$, denoted by $\dn(I)$, to be the simplicial complex over a 
set of vertices $\vs$, where $\{v_{i_1},\ldots, v_{i_s}\}$ is a face
of $\dn(I)$ if and only if $x_{i_1}\ldots x_{i_s} \notin I$. 
\end{itemize}
\end{definition}

\begin{notation}\label{notation} To simplify notation, we use $\DN$ to 
mean the nonface complex of $\F(\D)$ for a given simplicial complex $\D$. 
In other words, we set $$\DN = \dn(\F(\D)).$$ 

Given an ideal $I \subseteq k[V]$ where $V=\{\xs\}$, if there is no
reason for confusion, we use $\D$ and $\DN$ to denote $\df(I)$ and
$\dn(I)$, respectively.  If $F$ is a face of $\D=\tuple{\Fs}$, we let
the \emph{complements} of $F$ and $\D$ be
\begin{center}$\cmp{F}=V\setminus F$ and 
$\cmp{\D}=\tuple{\cmp{F_1},\ldots,\cmp{F_q}}$.
\end{center}
\end{notation}

\begin{definition}[Alexander dual]\label{alex} Let $I$ be a square-free
 monomial ideal in the polynomial ring $k[V]$ with $V=\{\xs\}$. Then
  the \emph{Alexander dual} of $\DN$ is the simplicial complex
  $$\alx{\DN} = \{ F \subset V \st \cmp{F} \notin \DN\}.$$
\end{definition}

It is easy to see that $\alx{\alx{\DN}} = \DN$.

We now focus on the relations between $\D$ and $\DN$ for a
given square-free monomial ideal $I$.  The first question we tackle is
how to construct $\DN$ from $\D$.

\begin{proposition}\label{duality-relations} Given a simplicial complex $\D$, 
  we have 
\begin{enumerate}
\item[(a)] $\DN=\cmp{\DM}$;

\item[(b)] $\alx{{\DN}}= {\DM}_N=\cmp{\D}$.
\end{enumerate}
\end{proposition}

\begin{proof} \begin{enumerate}

  \item[(a)] This is easy to check. See, for example, [BH] Theorem~5.1.4.
    
  \item[(b)] The last equality follows from  
             Proposition~\ref{cover-complex-dual} and
             Part (a), since $${\DM}{}_N = \cmp{{\DM}{}_M}=
             \cmp{\D}.$$
             
             We translate both sides of the first equation using the
             notations in~\ref{notation} and Definition~\ref{alex}:
             \begin{center} $\alx{\DN}=\{\cmp{F} \st F \notin \DN\}$ and
            $ \D^c= \tuple{\cmp{F} \st F {\rm \ is \ a \ facet\  of\ } \D}$.
            \end{center} 
             
            Suppose that $\cmp{F} \in \alx{\DN}$. Then $F \notin \DN$, and
            therefore if $f$ denotes the monomial that is the product
            of the vertices of $F$, and $I=\N(\DN)$, then $f \in I$.
            It follows that for some generator $g$ of $I$, $g |f$. If
            $G$ is the facet of $\D$ corresponding to $g$, we have $G
            \subseteq F$, which implies that $\cmp{F} \subseteq \cmp{G}$; so
            $\cmp{F} \in \cmp{\D}$.
            
            Conversely, let $G \in \cmp{\D}$. Then $G \subseteq
            \cmp{F}$, where $F$ is a facet of $\D$, so $f \in I$ which
            implies that $F \notin \DN$. So $\cmp{F} \in \alx{\DN}$,
            which implies that $G\in \alx{\DN}$.

          \end{enumerate}
        \end{proof}
           
        Proposition~\ref{duality-relations} is basically saying that
        the relationship between $\DM$ and $\alx{\DN}$ is the same as
        the relationship between $\D$ and $\DN$. The example below
        clarifies this point.

\begin{example}\label{diagram-example} Let $I=(xyz,zu) \subseteq k[x,y,z,u]$.
  Then the \emph{dual} ideal of $I$, which is the facet ideal of
  $\DM$, or equivalently the nonface ideal of $\alx{\DN}$, is the
  ideal $J=(xu,yu,z)$. The relationship between the four simplicial
  complexes and the two ideals is shown in Figure~\ref{diagram}.

\begin{figure}
\begin{center}
 \[ \includegraphics{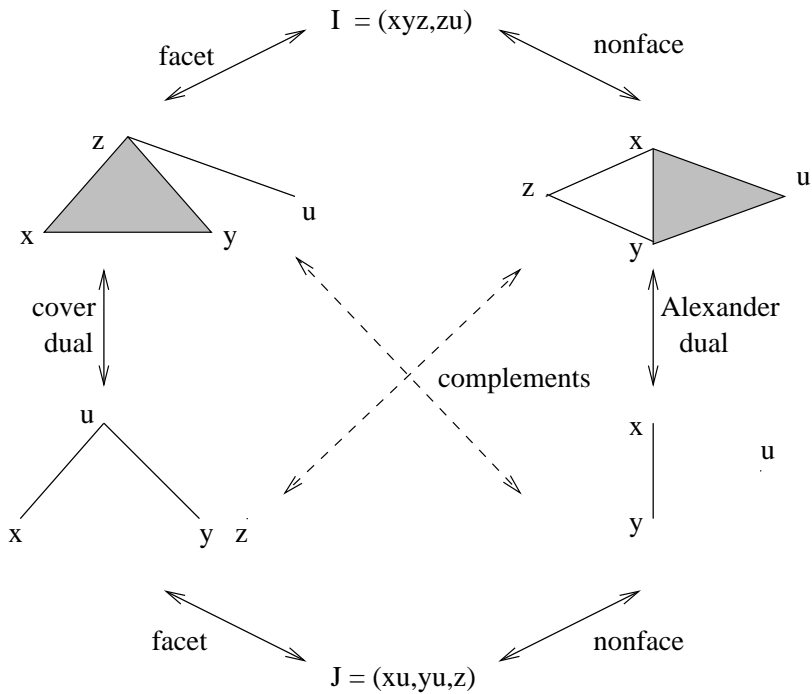} \]
\caption{Diagram of Example~\ref{diagram-example}}\label{diagram}
\end{center} 
\end{figure}
\end{example}

\idiot{Proposition~\ref{duality-relations} is basically just saying
what Ezra says in his combinatorics book-- that if $I=({\bf x}^{a_1},
\ldots, {\bf x}^{a_q})$ ($a_1,\ldots,a_q$ are integer vectors), and
$m$ is the maximal ideal $(\xs)$, then $\alx{I}=m^{a_1}\cap \ldots
\cap m^{a_q})$.}

Proposition~\ref{duality-relations} justifies the following definition.

\begin{definition}[dual of an ideal]\label{dual-ideal} Given a square-free 
  monomial ideal $I$ in a polynomial ring and $\D=\df(I)$, we define
  the $dual$ of $I$, denoted by $\IM$, to be the facet
  ideal of $\DM$, or equivalently, the nonface ideal of
  $\alx{\DN}$. So
  $$\IM=\F(\DM)=\N(\alx{\DN}).$$
\end{definition}

We now state a criterion for the Cohen-Macaulayness of a square-free
monomial ideal that is due to Eagon and Reiner ([ER]) in the language
stated above. First we define an ideal with a linear resolution.

\begin{definition}[linear resolution] An ideal $I$ in a polynomial ring 
  $R=k[\xs]$ over a field $k$, with the standard grading $\deg(x_i)=1$
  for all $i$, is said to have a \emph{linear resolution} if $R/I$ has
  a minimal free resolution such that for all $j>1$ the nonzero
  entries of the matrices of the maps $R^{\beta_j} \longrightarrow
  R^{\beta_{j-1}}$ are of degree 1.
\end{definition}

\begin{theorem}[{[ER]} Theorem 3]\label{eagon-reiner} Let $I$ be a 
square-free monomial ideal in a polynomial ring $R$. Then $R/I$ is
Cohen-Macaulay if and only if $\IM$ has a linear resolution.
\end{theorem}


\section{Simplicial Trees}\label{trees-section}

Considering simplicial complexes as higher dimensional graphs, one can
define the notion of a \emph{tree} by extending the same concept from
graph theory. Simplicial trees were first introduced in [F1] in order
to generalize results of [SVV] on facet ideals of graph-trees.  The
construction turned out to have interesting additional combinatorial
and algebraic properties.

Before we define a tree, we determine what ``removing a facet'' from a
simplicial complex means. We define this idea so that it corresponds
to dropping a generator from its facet ideal.

\begin{definition}[facet removal]\label{removal} Suppose $\Delta$
  is a simplicial complex with facets $\Fs$ and $\F(\Delta)=(\Ms)$ its
  facet ideal in $R=k[\xs]$. The simplicial complex obtained by
  \emph{removing the facet} $F_i$ from $\Delta$ is the simplicial
  complex
 $$\D \rmv{F_i}=\tuple{F_1,\ldots,\hat{F}_{i},\ldots,F_q}.$$
 \end{definition}
 
 Note that $\F(\D \rmv{F_i}) = (M_1,\ldots,\hat{M}_{i} ,\ldots,M_q)$.
 
 Also note that the vertex set of $\D \rmv{F_i}$ is a subset of the
 vertex set of $\D$.

\begin{example} Let $\D$ be a simplicial complex 
  with facets $F=\{x,y,z\}$, $G=\{y,z,u\}$ and $H=\{u,v\}$. Then $\D
  \rmv{F} =\tuple{G,H}$ is a simplicial complex with vertex set
  $\{y,z,u,v\}$.
\end{example}

In graph theory, a tree is defined as a connected cycle-free graph. An
equivalent definition is that a tree is a connected graph whose every
subgraph has a \emph{leaf}, where a leaf is a vertex that belongs to
only one edge. We make an analogous definition for simplicial
complexes by extending (and slightly changing) the definition of a
leaf.

\begin{definition}[leaf] A facet $F$ of a simplicial complex is called a 
  \emph{leaf} if either $F$ is the only facet of $\Delta$, or for some
  facet $G \in \Delta \rmv{F}$ we have
  $$F \cap (\D \rmv{F}) \subseteq G.$$
\end{definition}

Equivalently, the facet $F$ is a leaf of $\D$ if $F \cap (\D \rmv{F})$
is a face of $\D \rmv{F}$.

\begin{example}\label{leaf-example} Let $I=(xyz,yzu,zuv)$. Then $F=xyz$ is
 a leaf, but $H=yzu$ is not, as one can see in the picture below.

\[ \begin{tabular}{ccccc}
  \includegraphics{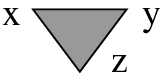} &$\cap$ &\includegraphics{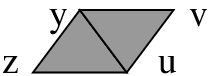} &= &
\includegraphics{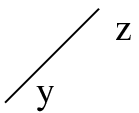} \\
&&&&\\
\includegraphics{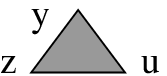} &$\cap$ &\includegraphics{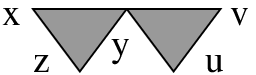} &= &
\includegraphics{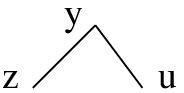}
\end{tabular} \]
\end{example}

\begin{definition}[tree, forest]\label{tree-def} A connected simplicial 
  complex $\Delta$ is a \emph{tree} if every nonempty subcollection of
  $\Delta$ has a leaf.  If $\D$ is not necessarily connected, but
  every subcollection has a leaf, then $\D$ is called a forest.
\end{definition}

\begin{example}\label{free-example} The simplicial complexes in 
examples~\ref{main-example} and~\ref{leaf-example} are both trees, but
  the one below is not because it has no leaves. It is an easy
  exercise to see that a leaf must contain a free vertex, where a
  vertex is \emph{free} if it belongs to only one facet.
\[ \includegraphics{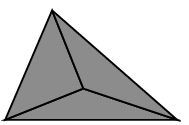} \]
\end{example}

An effective way to make algebraic arguments on trees is using
localization. It turns out that the minimal generating set of a
localization of the facet ideal of a tree corresponds to a forest.  As
we shall see below, this fact makes it easy to use induction on the
number of vertices of a tree. 

For details on the localization of s simplicial complex see [F2].
Here we give an example to clarify what we mean by localization.

\begin{example}\label{local-example} Let $\D$ be the simplicial complex below
  with $I=(xyz,yzu,yuv)$ its facet ideal in the polynomial ring
  $R=k[x,y,z,u,v]$.

\[ \includegraphics{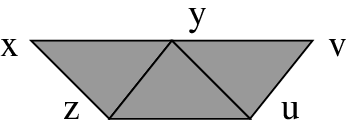} \]

Let $p=(x,u,z)$ be a prime ideal of $R$. Then $I_p=(xz,zu,u)=(xz,u)$
is the facet ideal of the forest below on the left. If $q= (y,z,v)$
then $I_q=(yz,yz,yv)=(yz,yv)$ corresponds to the tree on the right.

\[
\begin{tabular}{lclc}
Localization at $p$:&
\includegraphics{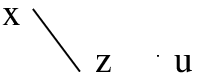}& 
\hspace{.2in}Localization at $q$:&
\includegraphics{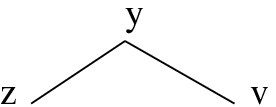} 
\end{tabular}
\]
\end{example}

Example~\ref{local-example} is an example of the following general
fact.

\begin{lemma}[Localization of a tree is a forest]\label{localization}
  Let $I \subseteq k[\xs]$ be the facet ideal of a tree, where $k$ is
  a field, and suppose that $p$ is a prime ideal of $k[\xs]$. Then for
  any prime ideal $p$ of $R$, $\df (I_p)$ is a forest. \end{lemma}

\begin{proof} See [F2] Lemma 4.5. 
\end{proof}


\section{Sequentially Cohen-Macaulay simplicial complexes}\label{sequential-section}

The notion of a sequentially Cohen-Macaulay ideal was introduced by
Stanley following the introduction of nonpure shellability by
Bj\"{o}rner and Wachs [BW]. It was known that every shellable
simplicial complex (which was by definition pure) was Cohen-Macaulay,
but what about nonpure shellable simplicial complexes? As it turns
out, ``sequentially Cohen-Macaulay'' is the correct notion to fill in
the gap here. On the other hand, the criterion of Eagon and Reiner
([ER]) stated that a simplicial complex is Cohen-Macaulay if and only if
its Alexander dual has a linear resolution. Herzog and Hibi ([HH])
developed the definition of a ``componentwise linear ideal'' so that
the above criterion extended to sequentially Cohen-Macaulay ideals: a
simplicial complex is sequentially Cohen-Macaulay if and only if its
Alexander dual is componentwise linear.

 In our setting, we use an equivalent characterization of sequentially
Cohen-Macaulay given by Duval, along with Theorem~\ref{eagon-reiner}
and the relationship between Alexander duality and cover complex
duality discussed in Section~\ref{dualities-section}, to prove that
simplicial trees are sequentially Cohen-Macaulay. In fact, we show
that if $I$ is the facet ideal of a simplicial tree, then the dual
$\IM$ of $I$ has ``square-free homogeneous components'' with
\emph{linear quotients}. This property is slightly stronger than what
we need, and it shows that if $I$ is a Cohen-Macaulay ideal to begin
with, then $\DN$ is shellable (which was known for the case where $\D$
is a graph; Theorem~6.4.7 of [V]).

Another outcome is the fact that an unmixed tree is Cohen-Macaulay
(Corollary~\ref{F2-cor}), which was shown in [F2] using very different
tools.

\begin{definition}[{[S]} Chapter III, Definition~2.9]\label{scm-def} 
  Let $M$ be a finitely generated ${\mathbb{Z}}$-graded module over a
  finitely generated ${\mathbb{N}}$-graded $k$-algebra, with $R_0=k$.
  We say that $M$ is \emph{sequentially Cohen-Macaulay} if there
  exists a finite filtration
  $$0=M_0 \subseteq M_1 \subseteq \ldots \subseteq M_r=M$$
  of $M$ by
  graded submodules $M_i$ satisfying the following two conditions.
\begin{enumerate}
\item[(a)] Each quotient $M_i/M_{i-1}$ is Cohen-Macaulay.

\item[(b)] $\dimn (M_1/M_0) < \dimn (M_2/M_1) < \ldots < \dimn (M_r/M_{r-1})$, 
where $\dimn$ denotes Krull dimension.
\end{enumerate}
\end{definition}

A simplicial complex is said to be \emph{sequentially Cohen-Macaulay}
if its Stanley-Reisner ideal has a sequentially Cohen-Macaulay
quotient.

The following characterization of a sequentially Cohen-Macaulay
simplicial complex given by Duval ([D] Theorem~3.3) is what we use in
this paper.

\begin{theorem}[{[D]} sequentially Cohen-Macaulay]\label{scm} Let $I$ be
  square-free monomial ideal $I$ in a polynomial ring $R$ over a field
  $k$, and let $\DN=\dn(I)$. Then $R/I$ is sequentially Cohen-Macaulay
  if and only if for every $i$, $-1 \leq i \leq \dimn \DN$, if $\DNi$
  is the pure $i$-dimensional subcomplex of $\DN$, then $R/\N(\DNi)$
  is Cohen-Macaulay.
\end{theorem}

\begin{example}\label{SCM-example} let $I=(xyz,zu)$ be the ideal 
of Example~\ref{diagram-example} in the diagram above.  Then for
 $i=0,1,2$, we have the following three simplicial complexes,
 respectively,
\[ \includegraphics{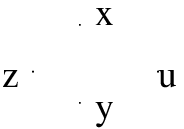} \hspace{1in}
\includegraphics{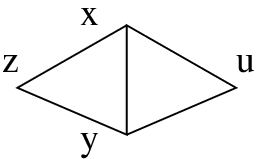} \hspace{1in}
\includegraphics{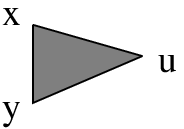} \]
which are, respectively, the nonface complexes of the ideals
$I_0=(xy,xz,xu,yz,yu,zu)$, $I_1=(xyz,xyu,zu)$ and $I_2=(z)$. One can
verify that all three of these ideals have Cohen-Macaulay quotients,
so $I$ is sequentially Cohen-Macaulay.
\end{example}

We define a componentwise linear ideal in the square-free case
using [HH]~Proposition~1.5.

\begin{definition}[square-free homogeneous component, 
componentwise linear]\label{comp-def} Let $I$ be a square-free
  monomial ideal in a polynomial ring $R$.  For a positive integer
  $k$, the $k$-th \emph{square-free homogeneous component} of $I$,
  denoted by $I_{[k]}$ is the ideal generated by all square-free
  monomials in $I$ of degree $k$.  The ideal $I$ above is said to be
  \emph{componentwise linear} if for all $k$, the square-free
  homogeneous component $I_{[k]}$ has a linear resolution.
\end{definition}

        Let $$\D=\tuple{\Fs}$$ be a simplicial complex with $\F(\D)
        \subseteq k[V]$, $V=\{\xs\}$, and let
        $$\DM=\tuple{\Gs}$$ be its cover complex . Then by
        Proposition~\ref{duality-relations} we know that 
        $$\DN=\tuple{\cmp{G_1},\ldots,\cmp{G_p}}.$$ For a given $i$, consider
        the pure $i$-dimensional subcomplex of $\DN$
        $$\DNi=\tuple{H_1,\ldots,H_u}.$$
        By Theorem~\ref{eagon-reiner}
        showing that $I_i= \N(\DNi)$ is a Cohen-Macaulay ideal is
        equivalent to showing that $\alx{I_i}$ has a linear resolution.
        By Proposition~\ref{duality-relations}, $\alx{I_i}$ is the facet
        ideal of $\cmp{\DNi}$.

        So we focus on $\cmp{H}$, where $H$ is a facet of
        $\DNi$. Since $H$ belongs to a subcomplex of $\DN$, for some
        facet $\cmp{G_j}$ of $\DN$, $H \subseteq \cmp{G_j}$. This
        implies that $G_j=\cmp{\cmp{G_j}} \subseteq \cmp{H}$;
        i.e. $\cmp{H}$ contains a minimal vertex cover of $\D$, and so
        $\cmp{H}$ is a variable cover of $\D$ of cardinality
        $n-(i+1)$.
        
        Similarly, if $G$ is a variable cover of cardinality
        $n-(i+1)$ of $\D$, then one can see that $\cmp{G}$ is a facet
        of $\DNi$.
        
        The discussion above shows that $\alx{I_i}$ is generated by
        monomials corresponding to variable covers of cardinality
        $n-i-1$ of $\D$. In other words
        $$\alx{I_i}={\alx{I}}_{[n-i-1]},$$ where ${\alx{I}}_{[j]}$
        denotes the $j$-th square-free homogeneous component of $\IM$,
        and showing that $\DNi$ is Cohen-Macaulay is equivalent to
        showing that ${\alx{I}}_{[n-i-1]}$ has a linear resolution.

        We have thus shown that:

\begin{proposition}[Criterion for being sequentially 
Cohen-Macaulay]\label{scm-criterion} Let $I$ be a square-free monomial
ideal in a polynomial ring. Then $I$ is a sequentially Cohen-Macaulay
ideal if and only if $\IM$ is componentwise linear.
\end{proposition}


\section{Simplicial trees are Sequentially Cohen-Macaulay}\label{trees-scm-section}

This section contains the main results of the paper. Our goal here is
to show that the facet ideal $I$ of a simplicial tree is sequentially
Cohen-Macaulay. By Proposition~\ref{scm-criterion} this is equivalent
to showing that the facet ideal $\IM$ of the cover complex of a tree
is componentwise linear (Definition~\ref{comp-def}). In fact, we show
that $\IM$ satisfies a stronger property: for every $i$, we show below
that $\IMi$ has linear quotients. This property, defined by Herzog and
Takayama in [HT], implies that $\IMi$ has a linear resolution. It also
implies additional combinatorial properties for $I$ (see
Corollary~\ref{shellabale-cor}).
 
\begin{definition}[linear quotients ({[HT]})]\label{linear-quotients-def} If 
$I \subset k[\xs]$ is a monomial ideal and $G(I)$ is its unique
minimal set of monomial generators, then $I$ is said to have
\emph{linear quotients} if there is an ordering $\Ms$ on the elements
of $G(I)$ such that for every $i=2,\ldots,q$, the quotient ideal
$$(M_1,\ldots,M_{i-1}):M_i$$ is generated by a subset of the variables
$\xs$.
\end{definition}

The following is a well-known fact. We reproduce an argument (almost
identical to one given in [Z] for the case of trees).

\begin{lemma}\label{quotient-resolution} If $I=(\Ms)$  is a monomial 
ideal in the polynomial ring $R=k[\xs]$ over the field $k$ that has
linear quotients and all the $M_i$ are of the same degree, then $I$
has a linear resolution.
\end{lemma}

              \begin{proof}  The proof is by induction on $q$. 
                The case $q=1$ is clear. Given that the ideal
                $I'=(M_1,\ldots,M_{q-1})$ has linear quotients and
                therefore a linear resolution, and that the degree of
                all the $M_i$ is equal to $d$, we have that (see
                Section~5.5 of [BH]) for all $i$:
                
                \begin{tabular}{ll} 
                $\tor{i}{R}(k,R/I')_a=0$ & unless $a= i+d$; \\
                $\tor{i}{R}(k,R/I':I)_a=0$ & unless $a=i+1$ ($I':I$ is
                generated by degree 1 monomials).
                \end{tabular}

                Consider the short exact sequence: $$0\longrightarrow
                R/(I':I) (-d) )\stackrel{.M_q}{\longrightarrow} R/I'
                \longrightarrow R/I \longrightarrow 0$$

               We obtain the long exact homology sequence:
               $$\cdots \longrightarrow \tor{i}{R}(k,R/(I':I)(-d))
               \longrightarrow \tor{i}{R}(k,R/I')\longrightarrow
               \tor{i}{R}(k,R/I)$$ $$ \longrightarrow
               \tor{i-1}{R}(k,R/(I':I)(-d)) \longrightarrow \cdots $$

               For a given $i$, $\tor{i}{R}(k,R/I)_a=0$ unless
               $$\tor{i}{R}(k,R/I')_a \neq0 \ {\rm or \ }
               \tor{i-1}{R}(k,R/(I':I))(-d)_a \neq 0.$$ Either way,
               this means that for any $i$, if $\tor{i}{R}(k,R/I)_a
               \neq 0$, then $a=i+d$. This implies that $R/I$ has a
               linear resolution.\end{proof}

We now set out to prove if $I \subseteq k[V]$, $V=\{\xs\}$, is the
  facet ideal of a tree (in fact, a forest) $\D$, and $i$, $\al(\D)
  \leq i \leq n$, is a given integer, then $\IMi$ has linear
  quotients.

We use induction on $n$. If $n=1$, $\D$ can only be the vertex
$\tuple{x_1}$, and so the only thing to check is if
${\IM}_{[1]}=(x_1)$ has linear quotients, which is obvious.

Suppose that $n>1$. We first deal with some special cases. If $\D$ is
a forest of singletons of the form $$\D=\tuple{x_1,\ldots,x_j}$$ where
$j<n$, then we can consider $I'=\F(\D)$ as an ideal in the polynomial
ring $R'=k[x_1,\ldots,x_{n-1}]$ ($I$ and $I'$ have the same
generating set, they only live in two different rings). By the
induction hypothesis, for every $i$, $\alx{I'} _ {[i]}$ has linear
quotients. 

It is easy to see that for every $i$, $$\IM_{[i]}= \alx{I'} _ {[i]} +
x_n \alx{I'} _ {[i-1]}.$$ Suppose that $$\alx{I'}_ {[i]}=(A_1,\ldots,
A_a) {\rm \ and \ } \alx{I'} _ {[i-1]}=(B_1,\ldots,B_b)$$ where the
generators of both ideals are written in the correct order for linear
quotients (recall that we are using the notation $xA$ to mean $\{x\}
\cup A$, since generally we are always thinking of sets as monomials).
To see that $$\IM_{[i]}=(A_1,\ldots, A_a) + x_n(B_1,\ldots,B_b)$$ has
linear quotients, we consider the case where for some monomial $m$ in
$k[\xs]$ (we can without loss of generality assume here that the products
are square-free),
$$mx_nB_j \in (A_1,\ldots, A_a,x_nB_1,\ldots,x_nB_{j-1}).$$ If
$mx_nB_j \in (x_nB_1,\ldots,x_nB_{j-1})$, since $\alx{I'} _ {[i-1]}$
has linear quotients, it follows that for some variable $z$ dividing
the monomial $m$, we have $zx_nB_j \in (x_nB_1,\ldots,x_nB_{j-1})$
(note that $m \neq 1)$.

If $mx_nB_j \in (A_1,\ldots, A_a)$, then since $B_j$ is already a
variable cover of $\D$, for any variable $z$ not in $B_j$, $zB_j$
covers $\D$ and is of cardinality $i$, and hence $zB_j \in
\{A_1,\ldots, A_a\}$. Therefore for any $z$ dividing $m$ we can again
conclude that $zx_nB_j \in (A_1,\ldots, A_a)$.

This argument settles the case where $\D=\tuple{x_1,\ldots,x_j}$, and $j<n$.

If $\D=\tuple{\xs}$, then the only ideal to consider is
$\IM_{[n]}=(x_1\ldots x_n)$ which by definition has linear quotients.
             
So now we can assume that $\D$ is a forest containing a facet with
more than one vertex.
    
 We begin our discussion with the following simple observation.
  
 \begin{lemma}\label{D'-lemma} Let $\D$ be  a simplicial complex with
$\F(\D)\subseteq k[V]$, $k$ a field, and $V=\{\xs\}$. Suppose that $x
   \in V$ is such that $V \setminus \{x\}$ is a variable cover for
   $\D$, and let $p_x$ be the prime ideal generated by the set $V
   \setminus \{x\}$. Then localizing $\D$ at $p_x$ corresponds, via
   the cover duality, to removing all facets of $\DM$ that contain
   $x$. In other words, if $\D'=\df(\F(\D)_{p_x})$ and $\As$ are the
   facets of $\DM$ that contain $x$, then $$\D'_M = \DM \rmv{\As}.$$
 \end{lemma}

               \begin{proof}  Note that a facet of
                 $\D'_M$ is the generating set for a minimal prime of
               $I=\F(\D)$ not containing $x$, and therefore belongs to
               $\DM$ as well. Conversely, if $A$ is a facet of the
               right-hand-side, then it corresponds to a minimal prime
               of $I$ not containing $x$ and hence to a minimal prime
               of $I_{p_x}$. \end{proof}

Now assume that the forest $\D$ has a leaf $F$ with positive dimension
and a free vertex (see Example~\ref{free-example}) $x=x_1$.  We can
write:
$$\DMi=\df(\IMi)= \tuple{\As} \cup \tuple{xB_1,\ldots,xB_s}$$ where
$\As$ are all the variable covers of $\D$ that have cardinality $i$
and do not contain $x$, and $xB_1,\ldots,xB_s$ are all the other
variable covers of cardinality $i$.

Now let $$\D'=\df(\F(\D)_{p_x}) \ {\rm \ and }\ \ \D''= \D \rmv{F}.$$
Both $\D'$ and $\D''$ are forests (by the definition of a tree, and by
Lemma~\ref{localization}) whose vertex sets are contained in
$\{x_2,\ldots,x_n\}$. Also note that $\D'$ is a nonempty simplicial
complex.

With notation as above, by Lemma~\ref{D'-lemma} $$\D'_{M,\
[i]}=\tuple{\As}.$$

Also notice that $$\D''_{M,\ [i-1]}=\tuple{\Bs}.$$ To see this last
equation, note that since for $j=1,\ldots,s$, $xB_j$ covers $\D$,
$B_j$ has to cover $\D''$ (as $x$ is a free vertex of $F$ and hence
only covers $F$). On the other hand, if $A$ is any variable cover of
$\D''$ of cardinality $i-1$, then $xA$ is in $\DMi$, and so $xA \in
\{xB_1,\ldots,xB_s\}$.

Applying the induction hypothesis to the forests $\D'$ and $\D''$ we
see that the ideals
\begin{center} 
${\alx{I'}}_{[i]}=(\As)$ and ${\alx{I''}}_{[i-1]}=(\Bs)$ 
\end{center}
of $k[x_2,\ldots,x_n]$ both have linear quotients. Without loss of
generality assume that the given orders on the $A$'s and the $B$'s are
appropriate for taking quotients.  We show that the ideal
$${\IM}_{[i]}=(\As)+ x(\Bs)$$ also has linear quotients.  Here we
assume that $1 < i <n$, since ${\IM}_{[n]}=(x_1 \ldots x_n)$ has
linear quotients by definition, as does ${\IM}_{[1]}$ which is, if
nonzero, generated by a subset of $\{\xs\}$.

The first case of interest is the ideal $$(\As): xB_1.$$ Now $B_1$ is
a variable cover of $\D''= \D \rmv{F}$, so $yB_1 \in {\alx{I}}_{[i]}$
for any vertex $y$ of $F$ not in $B_1$. So if $m$ is any monomial such
that $mxB_1 \in {\alx{I'}}_{[i]}$, then for some monomial $n$ and some
$j$, assuming without loss of generality that both products below are
square-free, we have
$$mxB_1 =nA_j.$$ If $B_1$ already contains a vertex of $F$, then it is
a variable cover of cardinality $i-1$ for $\D'$, and so for any $y|m$,
$yB_1 \in \{\As\}$. Otherwise, since there is some vertex $y$ of $F$
in $A_j$, $y$ has to divide $m$, which again implies that $yxB_1 \in
(\As)$.

In general, for the ideal
$$(\As,xB_1,\ldots,xB_{j-1}):xB_j$$ if for some monomial $m$, $mxB_j
\in (xB_1,\ldots,xB_{j-1})$, then by the induction hypothesis on
${\alx{I''}}_{[i-1]}$ there is a variable $y$ that divides $m$ such that
$yxB_j \in (xB_1,\ldots,xB_{j-1})$.

If $mxB_j \in (\As)$, then it follows from an argument identical to
the case $j=1$ above that there is a variable $y$ dividing $m$ such that
$yxB_j \in (\As)$.

\idiot{ If $\D$ has no variable cover of order $i$, we assume
$\alx{I}_{[i]}=(0)$ which has linear quotients as it is a principal
ideal.}

We have thus proved that:

\begin{theorem}\label{linear-quotient-theorem} If $I \subseteq k[\xs]$ is
the facet ideal of a simplicial tree (forest) $\D$, then $\IMi$ has
linear quotients for all $i=\al(\D),\ldots,n$.
\end{theorem}

Theorem~\ref{linear-quotient-theorem} along with
Lemma~\ref{quotient-resolution} result in the following statement.

\begin{corollary}\label{componentwise-corollary} If $\D$ is a simplicial 
tree (forest), then $\alx{\F(\D)}$ is a componentwise linear
 ideal. \end{corollary}

Putting Corollary~\ref{componentwise-corollary} together with
Proposition~\ref{scm-criterion}, we arrive at our final goal.

\begin{corollary}[Trees are sequentially Cohen-Macaulay]\label{main-result}
  The facet ideal of a simplicial tree (forest) is sequentially
  Cohen-Macaulay.
\end{corollary}

\begin{example} The ideal $I$ in Example~\ref{SCM-example} is sequentially  
Cohen-Macaulay because it is the facet ideal of a tree. 
\end{example}

It follows easily that if the tree $\D$ is unmixed to begin with, then
it must be Cohen-Macaulay. This is because in this case $\alx{\F(\D)}$
itself is a square-free homogeneous component, which has a linear
resolution. So by applying Theorem~\ref{eagon-reiner} we have

\begin{corollary}[An unmixed tree is Cohen-Macaulay]\label{F2-cor}
  If $\D$ is an unmixed simplicial tree, then $\F(\D)$ has a
  Cohen-Macaulay quotient.  \end{corollary}

Corollary~\ref{F2-cor} was proved in [F2] using very different
tools. In particular, in [F2] we show that a tree is unmixed if and
only if it is ``grafted''. The notion of grafting is what gives a
Cohen-Macaulay tree its definitive combinatorial structure.

Another interesting fact that follows is that in the case of a
simplicial tree $\D$, if $\D$ is Cohen-Macaulay, then $\DN$ is
shellable (see [BH] for the definition). Given a square-free monomial
ideal $I$, if $\dn(I)$ is shellable, then $I$ is Cohen-Macaulay (see
[BH]), but the converse is not true in general.

\begin{corollary}\label{shellabale-cor} If 
  $\D$ is a Cohen-Macaulay simplicial tree, then $\D_N$ is shellable.
\end{corollary}

       \begin{proof} If $I=\F(\D)$ is Cohen-Macaulay, then by 
        Theorem~\ref{linear-quotient-theorem}, $\IM$ has linear
        quotients (since it has generators of the same degree). The
        rest follows directly from the definitions of shellability and
        linear quotients; see [HHZ]~Theorem~1.4, part~(c). \end{proof}

\begin{figure}

\begin{tabular}{|l|l|l|}
\hline &&\\ 
{\bf Notation} & {\bf Meaning} & {\bf First appearance}\\
&&\\ 
\hline &&\\
 $\F(\D)$& facet ideal of $\D$ & Definition~\ref{facet-ideal}\\ 
$\df(I)$& facet complex of $I$& Definition~\ref{facet-ideal}\\ 
$\al(\D)$ & vertex covering number of $\D$ & Definition~\ref{unmixed}\\
$\DM$ & cover complex of $\D$ & Definition~\ref{cover-complex}\\ 
$\N(\D)$& nonface ideal of $\D$ & Definition~\ref{nonface-ideal}\\ 
$\dn(I)$& nonface complex of $I$& Definition~\ref{nonface-ideal}\\ 
$\DN$ & $\dn(\F(\D))$ & Notation~\ref{notation} \\ 
$\cmp{F}$, $\cmp{\D}$ & complements of $F$ and $\D$ & 
Notation~\ref{notation}\\ 
$\alx{\D}$ & Alexander dual of $\D$ & Definition~\ref{alex}\\ 
$\IM$& dual of $I$ & Definition~\ref{dual-ideal}\\ 
$\D \rmv{F}$ & removal of facet $F$ from $\D$ & Definition~\ref{removal}\\ 
$\DNi$ & pure $i$-dimensional subcomplex of $\DN$ & Theorem~\ref{scm}\\ 
$I_{[i]}$& $i$-th square-free homogeneous component of $I$ & 
Definition~\ref{comp-def}\\
$p_x$ & ideal generated by all variables but $x$ &
Lemma~\ref{D'-lemma}\\ 
$\DMi$ & facet complex of $\IMi$ & following Lemma~\ref{D'-lemma}\\ 
&&\\ \hline
\end{tabular}
\caption{Index of Notation}\label{index}
\end{figure}



\begin{thebibliography}{1234}
  
\bibitem[B]{B} Berge, C. \emph{Hypergraphs, Combinatorics of finite
    sets}, North-Holland Mathematical Library, 45. North-Holland
  Publishing Co., Amsterdam, 1989.
  
\bibitem[BH]{BH} Bruns, W., Herzog, J. \emph{ Cohen-Macaulay
rings}, vol.  39, Cambridge studies in advanced mathematics, revised
edition, 1998.

\bibitem[BW]{BW} Bj\"{o}rner, A., Wachs, M.L. \emph{Shellable nonpure
    complexes and posets}, I.  Trans. Amer. Math. Soc.  348 (1996),
  no. 4, 1299--1327.

\bibitem[D]{D} Duval, A.M. \emph{Algebraic shifting and 
sequentially Cohen-Macaulay simplicial complexes}, Electron. J.
  Combin. 3  (1996), no. 1, Research Paper 21

\bibitem[ER]{ER} Eagon J.A., Reiner, V. \emph{Resolution of
Stanley-Reisner rings and Alexander duality}, J. Pure Appl. Algebra
130 (1998), no. 3, 265--275.

\bibitem[F1]{F1} Faridi, S. \emph{ The facet ideal of a simplicial
complex}, Manuscripta Mathematica 109 (2002), 159-174.

\bibitem[F2]{F2} Faridi, S. \emph{Cohen-Macaulay properties of
    square-free monomial ideals}, Preprint.
  
\bibitem[HH]{HH} Herzog, J., Hibi, T. \emph{Componentwise
    linear ideals}, Nagoya Math. J. 153 (1999), 141--153.  

\bibitem[HHZ]{HHZ} Herzog, J., Hibi, T., Zheng, X. \emph{Dirac's
theorem on chordal graphs and Alexander duality}, Preprint.

\bibitem[HRW]{HRW} Herzog, J., Reiner, V., Welker,
V. \emph{Componentwise linear ideals and Golod rings},  Michigan
Math. J.  46 (1999), no. 2, 211--223.


\bibitem[HT]{HT} Herzog, J., Takayama, Y. \emph{Resolutions by mapping
    cones}, The Roos Festschrift volume, 2.  Homology Homotopy Appl. 4
  (2002), no. 2, part 2, 277--294 (electronic).

\bibitem[S]{S} Stanley, R.P. \emph{Combinatorics and commutative
    algebra}, Second edition. Progress in Mathematics, 41. Birkhäuser
  Boston, Inc., Boston, MA, 1996. x+164 pp. ISBN: 0-8176-3836-9.
  
\bibitem[SVV]{SVV} Simis A., Vasconcelos W.,
  Villarreal R., \emph{On the ideal theory of graphs}, J. Algebra 167
  (1994), no. 2, 389--416.

\bibitem[V]{V} Villarreal R., \emph{Monomial algebras},
Monographs and Textbooks in Pure and Applied Mathematics, 238. Marcel
Dekker, Inc., New York, 2001.

\bibitem[Z]{Z} Zheng, X. \emph{Resolutions of facet ideals}, Preprint.


\end{thebibliography}
\end{document}